\documentclass[runningheads]{llncs}
\usepackage[utf8]{inputenc}
\usepackage{amsmath,amssymb}

\usepackage{enumitem}
\usepackage{algorithm}
\usepackage[noend]{algpseudocode}
\usepackage[hidelinks]{hyperref}

\makeatletter
\renewcommand{\endproof}{\qed\@endtheorem}
\makeatother

\newenvironment{namedproof}[1]{\begin{proof}[#1]}{\end{proof}}

\newcommand{\R}{\mathbb R}

\newcommand{\zero}{\mathbf 0}
\newcommand{\supp}{\operatorname{supp}}
\newcommand{\T}{^{\mathsf T}}

\newcommand{\NN}{N}
\newcommand{\Cost}[2]{\mathcal C(#1,#2)}
\setlist[enumerate]{leftmargin=*,itemsep=2pt,topsep=4pt}
\allowdisplaybreaks[1]

\begin{document}

\title{Polynomial sequences of Simplex pivots\\on $0/1$-polytopes}
\titlerunning{Polynomial sequences of Simplex pivots on $0/1$-polytopes}
\author{Laura Sanit\`a} 
\authorrunning{Laura Sanit\`a}
\institute{Bocconi University of Milan (Italy)\\
\email{laura.sanita@unibocconi.it}}
\maketitle

\begin{abstract}
We study the number of Simplex pivots needed to solve linear
programs over $0/1$-polytopes, starting from an arbitrary feasible basis
and counting degenerate pivots. It has been known since the 1990s that from any extreme point of a
$0/1$-polytope there exists a monotone path of linear length to an optimal one. However, translating this into a polynomial
sequence of Simplex pivots, i.e., improving basis exchanges, has remained an open question.

We here provide an affirmative answer, showing that for every $0/1$-polytope $\mathcal P \subseteq \mathbb R^d$, an optimal
basis can be reached with $O(d^3)$ Simplex pivots, 
independently of the number of inequalities and of the size of the
coefficients.

\end{abstract}

\section{Introduction}\label{sec:intro}

The Simplex method, introduced by Dantzig~\cite{Dantzig}, is among the
most widely used algorithms for solving linear programs (LPs). Its
influence on science and engineering led to its inclusion among the
ten most influential algorithms of the twentieth century~\cite{DS}.
Besides its practical importance, the Simplex method remains a
candidate for resolving one of the central open questions in
discrete mathematics: whether linear programming admits a strongly polynomial
time algorithm~\cite{Smale}. However, despite decades of research,
its theoretical performance is still poorly understood, even for classes of LPs with considerable
combinatorial structure.

The Simplex method starts at a \emph{feasible basis} of an LP, which corresponds to a
\emph{minimal set of tight inequalities} representing an
extreme point solution. At each iteration, it exchanges bases by swapping one
 inequality with a different one: this operation geometrically corresponds to selecting a direction $y$ that is
\emph{improving} according to the objective function.
Such an exchange is called a \emph{Simplex pivot}. If a positive step can be taken in the selected direction,
the algorithm moves along an edge of the feasible region to a better
extreme point. Otherwise, the basis changes while the extreme point
remains the same, and the pivot is called \emph{degenerate}. Thus, the
sequence of distinct extreme points visited by the algorithm traces a
\emph{monotone} path on the $1$-skeleton\footnote{The $1$-skeleton of a polyhedron is the graph whose nodes are its extreme points and whose edges are its one-dimensional faces. A path in this graph is monotone if the objective value strictly increases along each extreme point of the path.} of the feasible polyhedron, but because of degeneracy, the
number of pivots may be much larger than the length of this path.

This distinction is central to the complexity of the Simplex method.
Its performance is closely related to the so-called \emph{polynomial
Hirsch conjecture}, which asks whether the shortest path between any
two extreme points on the $1$-skeleton of a polyhedron can be bounded
by a polynomial in its dimension and number of facets~\cite{Santos}.
This is a major open question for general polyhedra, but for several
important classes polynomial diameter bounds are known, and one can
even construct monotone paths of polynomial length from any initial
vertex to an optimum. Degeneracy is often avoided by perturbing the
LP, but a perturbation may destroy precisely the structure responsible
for such bounds. Thus, whether this geometric information can be
translated into a polynomial number of Simplex pivots, that is, of
improving basis exchanges on the original formulation, is in general
not clear.

A further subtlety matters for our paper. At a degenerate vertex
one could allow arbitrary degenerate basis exchanges, using them only to reach a basis from which the
desired edge can be followed: with such exchanges, translating a short
monotone path into a short sequence of basis changes is trivial. 
We insist instead that every exchange be a \emph{Simplex
pivot}: the entering variable must have positive reduced cost, even
when the pivot is degenerate and does not move the current solution.
This is essential if the goal is to understand the Simplex method
itself, and it is the analogue, at the level of bases, of requiring
paths on the $1$-skeleton to be monotone: every step must be justified
by the objective function. The question we ask is therefore
deliberately stronger: given a short monotone path, can its progress
be reproduced by a polynomial sequence of Simplex pivots, degenerate
ones included?

A prominent example is the class of $0/1$-LPs, namely LPs defined on $0/1$-polytopes, 
whose vertices have all their coordinates in $\{0,1\}$. Such polytopes
arise throughout combinatorial optimization, and their geometry is
particularly favorable. The Hirsch bound holds for them, even with respect to monotone paths~\cite{Naddef,MatsuiTamura},
and several geometric rules yield short monotone paths, including
greatest-improvement~\cite{Chubanov}, steepest-edge and shadow rules~\cite{DLS,BDLKS,BKS}.
In other words, for $0/1$-LPs, the existence of a monotone path of strongly polynomial length from
any arbitrary initial feasible extreme point to an optimal one has been known since the 90s.
However, the existence of a sequence of
Simplex pivots of strongly polynomial length from an arbitrary initial feasible basis to an optimal one has remained
open for all these years. This shows a clear gap in our understanding of the Simplex method even for
a very classical setting, and in fact the gap between short edge paths and short pivot sequences for $0/1$-LPs has been
explicitly mentioned as an open question in the literature.

\subsection{Our results and comparison with the literature}
In this paper, we finally give an affirmative answer.  Our main result is the following.

\begin{theorem}\label{thm:main}
Let $\max\{c\T x: x\in\mathcal P\}$ be a $0/1$-LP with
$\mathcal P\subseteq\R^d$. From every feasible initial basis, an
optimal basis can be reached by a sequence of at most $O(d^3)$
Simplex pivots, at most $O(d^2)$ of which are nondegenerate.
\end{theorem}

The bound is independent of the number of inequalities describing
$\mathcal P$ and of the size of the coefficients: our
arguments are purely combinatorial. In
particular, for rational data, the bound does not depend on the
encoding size. We remark that here
$\mathcal P\subseteq\R^d$ is a $0/1$-polytope given by any system of
linear \emph{inequalities}.  
The Simplex method can operate on it after putting it in Standard
Equality Form (SEF) (the precise definitions are given in
Section~\ref{sec:prelim}). Assuming that in
SEF the description has $n$ variables and $m$ equality constraints, where
$n-m\leq d$, the exact bound we prove is $\sum_{r=1}^{n-m}r\min\{r,m\}+n-m$
(see the proof of Theorem~\ref{thm:main} in
Section~\ref{sec:algorithm}).

Our starting point is the antistalling procedure of~\cite{KS}, which
escapes a degenerate vertex of an LP in standard equality form with $n$
variables and $m$ constraints in at most $n-m-1$ degenerate Simplex
pivots, given an improving feasible direction $y$ as a guide. This does
not translate into a bound on the number of Simplex pivots, even
when $y$ is the direction of an edge of a short monotone path whose
existence is guaranteed: the procedure may modify $y$ during its
execution, and in any case it does not guarantee that the
nondegenerate pivot it eventually performs follows $y$. Nevertheless,
\cite{KS} obtains a polynomial bound on the  number of Simplex
pivots for \emph{some} $0/1$-LPs (namely, LPs that once put in
SEF have nonzero coordinates whose ratio is polynomially bounded), by showing that
every nondegenerate pivot makes substantial progress in the objective
function, in the spirit of the work of~\cite{KM}. The extension to all $0/1$-LPs is left open.

The authors of~\cite{Chubanov,BKS,DLS,BDLKS} studied greatest-improvement, steepest-edge and shadow-vertex rules, whose monotone paths
are proven to be short on $0/1$-polytopes. In particular, the authors of \cite{BKS} show
that a prescribed shadow edge can be followed after linearly many
degenerate pivots, provided that the initial basis satisfies a
compatibility condition, but they do not obtain a bound on the total
number of Simplex pivots, which remains open. In fact, they show
that following a given edge may be impossible: they exhibit a
$0/1$-LP and a feasible basis at which the only Simplex
pivot is nondegenerate and traverses an edge with an arbitrarily small objective function gain compared to the
desired one. 

We here overcome this difficulty by not measuring progress through the objective improvement of one individual
nondegenerate pivot, but rather measuring progress through the size of common faces. 
We fix a target vertex $t$ (an optimal vertex for
Theorem~\ref{thm:main}, but any better vertex suffices) and always use
$t-x$ as the guide at the current vertex $x$. When a nondegenerate
pivot makes a variable zero that is also zero at $t$, we choose it as
the leaving variable and \emph{lock} it, that is, we never let it enter
again; this confines the subsequent pivots to a smaller face that still
contains $t$. When this is impossible, we show that the smallest face
$F(x,t)$ containing $x$ and $t$ strictly \emph{grows}, using a simple property
of $0/1$-polytopes: a vertex $x$ is determined by the face $F(x,t)$
(Lemma~\ref{lem:binary}).  Together with
the bound of~\cite{KS}, this gives a sequence of at
most $\sum_{r=1}^{n-m}r\min\{r,m\}$ improving pivots after which the
objective value of $t$ is attained or exceeded
(Theorem~\ref{thm:target}).

Since any better vertex can serve as a target, Theorem~\ref{thm:target}
also shows that the geometric information provided by short monotone
paths can always be exploited: the objective function gain of any
monotone path of length $L$ can be simulated by Simplex pivots with an
overhead of $O(d^3)$ pivots per edge (see Corollary~\ref{thm:simulate}). We
stress that the path is not followed: its vertices are only used as
targets, one after the other, and the vertices actually visited by the
pivots are in general different. This has an important consequence: namely,
the path need not be known in advance, nor obtained by a global
computation. When a \emph{local} rule for generating the path is available,
it can be computed vertex by vertex, using at each of its vertices
only local information such as the constraints tight there,
even though the Simplex method need not visit these vertices. In
other words, Theorem~\ref{thm:target} can be used to turn any rule that selects
the next vertex of a short monotone path from the local geometry at
the current vertex into a Simplex pivot rule, with a polynomial overhead.

The key point in all of this is that, in order to recover the
gain of the target, we must allow individual pivots that look useless:
some may improve the objective by a negligible amount, others seem to move
away from the target, as the smallest face containing both grows, instead of shrinking.
However, a growing face is in fact still a useful  measure of progress,
and the improvement of the target is anyway recovered after a cubic number of Simplex pivots.

\section{Preliminaries}\label{sec:prelim}

\paragraph{Linear programs and Simplex pivots.}
For an integer $k$, we let $[k]=\{1,\ldots,k\}$. For a vector $v$, we let $\supp(v)=\{i: v_i\neq0\}$,
and $\zero$ denotes the zero vector. For a matrix $A$, a vector $x$, and a subset of indices $B$, we let $A_B$ and $x_B$
denote the submatrix of $A$ and the subvector of $x$ indexed by $B$.
We consider linear programs in Standard Equality Form to be of the form
\begin{equation}\label{eq:LP}
 \max\{c\T x: Ax=b,\ x\geq 0\},
\end{equation}
where $A\in\R^{m\times n}$ has full row rank $m$, $b\in\R^m$, $c\in\R^n$,
and the feasible region $Q=\{x\in\R^n: Ax=b,\ x\geq0\}$ is a nonempty
polytope. A \emph{basis} is a set $B\subseteq[n]$ of $m$ indices with
$A_B$ nonsingular, and $\NN=[n]\setminus B$ is the set of \emph{nonbasic}
indices. The \emph{basic solution} associated to $B$ is given by $x_B=A_B^{-1}b$
and $x_{\NN}=\zero$; $B$ is \emph{feasible} if $x_B\geq\zero$, and $x$
is then a \emph{basic feasible solution} (BFS), which corresponds to a vertex of $Q$.
A feasible basis is \emph{degenerate} if $x_i=0$ for some $i\in B$;
a degenerate vertex may be represented by several feasible bases. For a
feasible basis $B$ and $j\in\NN$, the \emph{basic direction}
$r^j\in\ker(A):=\{y \in \mathbb R^n: Ay =0\}$ is defined by $r^j_j=1$, $r^j_i=0$ for
$i\in\NN\setminus\{j\}$, $r^j_B=-A_B^{-1}A_j$, and the \emph{reduced
cost} of $j$ is $\bar c_j=c\T r^j$. A \emph{pivot} with
entering index $j$ moves from $x$ along $r^j$ with step length
\begin{equation}\label{eq:ratio}
 \theta=\min\left\{x_i/(-r^j_i): i\in B,\ r^j_i<0\right\},
\end{equation}
where the set is nonempty as $Q$ is bounded; any index $g$ attaining
the minimum can be chosen as \emph{leaving} index, and the new basis
$(B\setminus\{g\})\cup\{j\}$ is feasible with BFS $x'=x+\theta r^j$.
The pivot is \emph{degenerate} if $\theta=0$ and \emph{nondegenerate}
otherwise. A \emph{Simplex pivot} is a pivot with $\bar c_j>0$ and with a leaving
index attaining the minimum in~\eqref{eq:ratio}. 

A direction $y\in\ker(A)$ is \emph{improving} if $c\T y>0$ and
\emph{feasible} at $x\in Q$ if $y_i\geq0$ whenever $x_i=0$. For a
feasible basis $B$, every $y\in\ker(A)$ satisfies
$y=\sum_{j\in\NN}y_jr^j$, hence
\begin{equation}\label{eq:rcidentity}
 c\T y=\textstyle\sum_{j\in\NN}\bar c_j\,y_j ;
\end{equation}
in particular, if $y$ is improving and feasible at the BFS of $B$, some
$j\in\NN$ has $y_j>0$ and $\bar c_j>0$. A feasible basis is
\emph{optimal} if $\bar c_j\leq0$ for all $j\in\NN$.

\paragraph{Faces.}
A face of a polytope $P$ is the empty set, $P$ itself, or any subset formed by the intersection of the polytope with a supporting hyperplane. 
For points
$x,y$ of a polytope, $F(x,y)$ denotes the smallest face containing
both. For $x,y\in Q$,
$Z(x,y)=\{i\in[n]: x_i=y_i=0\}$ is the set of their \emph{common
zeros}.

\paragraph{$0/1$-LPs.}
A \emph{$0/1$-polytope} is a polytope $\mathcal P\subseteq\R^d$ with all
vertices in $\{0,1\}^d$. We consider a \emph{$0/1$-LP} to be
a problem of the form $\max\{c\T v: v\in \mathcal P\}$ with $\mathcal P$ a $0/1$-polytope given by a general
system
\begin{equation}\label{eq:P}
 \mathcal P=\{v\in\R^d: Ev=f,\ Hv\leq h,\ v\geq\zero\},
  \qquad E\in\R^{e\times d},\ H\in\R^{\ell\times d},
\end{equation}
where $E$ has full row rank $e$ (possibly $e=0$). The nonnegativity constraints are without loss
of generality, as $v\geq\zero$ is valid for every $0/1$-polytope. The
system may contain redundant inequalities. The $0/1$-LP put in Standard Equality Form becomes
\begin{equation}\label{eq:SEF}
 \max\{c\T v: Ev=f,\ Hv+s=h,\ v\geq\zero,\ s\geq\zero\},
\end{equation}
an LP of the form~\eqref{eq:LP} in the variables $x=(v,s)$, with
$n=d+\ell$, $m=e+\ell$, $A=\left(\begin{smallmatrix}E&0\\H&I_\ell\end{smallmatrix}\right)$,
$b=(f,h)$, and objective $(c,\zero)$, which we also denote by $c$;
note that $n-m=d-e\leq d$. We denote by $V$ the index set of the
original variables, so that $x_V=v$. The map
$x\mapsto x_V$ is an affine bijection between the feasible region $Q$
of~\eqref{eq:SEF} and $\mathcal P$, with inverse $v\mapsto(v,h-Hv)$, which
preserves faces and vertices; in particular, every vertex $x$
of $Q$ is determined by its binary part $x_V\in\{0,1\}^d$, while its
slack coordinates need not be binary.

\paragraph{The antistalling procedure of~\cite{KS}.}
Given a feasible basis with BFS $x$ and an improving direction $y$
feasible at $x$ (the \emph{guide}), the antistalling procedure
of~\cite{KS} selects an entering index $f\in\NN$ with $y_f>0$ and
$\bar c_f>0$, which exists by~\eqref{eq:rcidentity}, performs a
few degenerate pivots, after replacing $y$ by $y+\alpha r^f$ for
a suitable $\alpha>0$ if needed, and repeats with the new basis until
a nondegenerate pivot occurs. Each step
uses a polynomial number of arithmetic operations. We use the
following escape bound, stated for the maximization problem~\eqref{eq:LP}.

\begin{lemma}[{\cite[Theorem~2]{KS}}]\label{lem:KSescape}
Let $B$ be a feasible basis of~\eqref{eq:LP} with BFS $x$, and let $t$
be a vertex of $Q$ with $c\T t>c\T x$. The antistalling procedure with
guide $t-x$ performs a nondegenerate Simplex pivot after at most 
$\min\{n-m-1,m-1\}$ degenerate
ones.
\end{lemma}

Lemma~\ref{lem:KSescape} only concerns the pivots up to the first
nondegenerate one, hence at that pivot any index attaining the minimum
in~\eqref{eq:ratio} can be chosen as the leaving index.

\section{Attaining the value of a target vertex}\label{sec:target}

Throughout this section we consider a $0/1$-LP that has been put in SEF~\eqref{eq:SEF}. Let $B^0$ be a feasible basis with BFS $x^0$, and
let $t$ be a vertex of $Q$ with $c\T t>c\T x^0$, the \emph{target}.

\begin{theorem}\label{thm:target}
Let $\Cost{n}{m}=\sum_{r=1}^{n-m}r\min\{r,m\}$. From $B^0$ one can
reach a feasible basis corresponding to a BFS $x$ which satisfies $c\T x\geq c\T t$,
with a sequence of at most $\Cost{n}{m}$ Simplex pivots. The number of nondegenerate ones
is at most $(n-m)(n-m+1)/2$.
\end{theorem}

We remark that the vertex reached need not be $t$. The sequence is produced by the
antistalling procedure with guide $t-x$, re-initialized after every
nondegenerate pivot. Since the objective improvement of a nondegenerate
pivot can be negligible, as discussed in the introduction, we measure
progress through two faces of $Q$ containing $t$: a face obtained by
\emph{locking} variables, and the face $F(x,t)$. 

\paragraph{Locking variables.}
Let $B$ be a feasible basis with BFS $x$, and let $J\subseteq[n]$ be a
set of \emph{locked} indices such that
\begin{equation}\label{eq:J}
 J\cap B=\varnothing\qquad\text{and}\qquad x_j=t_j=0\quad(j\in J).
\end{equation}
Let $Q_J=\{z\in Q: z_j=0\ (j\in J)\}$, a face of $Q$ containing $x$ and
$t$, and consider the \emph{restricted LP} obtained from~\eqref{eq:LP}
by deleting the columns indexed by $J$, whose feasible region is $Q_J$.
Since $A_{[n]\setminus J}$ contains the nonsingular submatrix $A_B$,
the restricted LP is again of the form~\eqref{eq:LP}, with $n-|J|$
variables and $m$ constraints. We let $r(J)=n-|J|-m$ and note that
the dimension of $ Q_J$ is $\leq r(J)$. The following proposition is immediate.

\begin{proposition}\label{lem:restrict}
Let $B$ be a feasible basis of~\eqref{eq:LP} and let $J$
satisfy~\eqref{eq:J}. Then 
\begin{itemize}
\item[(i)] $B$ is a feasible basis of the restricted
LP with the same BFS $x$;
\item[(ii)] Every Simplex pivot of the restricted LP is a Simplex
pivot of~\eqref{eq:LP} whose entering index is not in $J$;
\item[(iii)] $t$ is a vertex of $Q_J$, and hence a BFS of the restricted LP.
\end{itemize}
\end{proposition}

\begin{proof}
The matrix $A_B$, the vector $A_B^{-1}b$ and the column $A_j$ do
not change when the columns in $J$ are deleted; hence $B$ is a feasible
basis of the restricted LP with BFS $x$, proving $(i)$. 
For $j\notin B\cup J$ the
direction $r^j$ and $\bar c_j=c\T r^j$ are
unchanged. The ratio test~\eqref{eq:ratio} only involves the basic
entries of $x$ and $r^j$, which are unchanged as well. This yields $(ii)$. 
Finally, $t\in Q_J$ because $t_J=\zero$, and a vertex of $Q$ lying in a face of
$Q$ is a vertex of that face. This proves $(iii)$.
\end{proof}

Thus pivoting in the restricted LP simply means refusing the locked
indices as entering indices. Moreover, Lemma~\ref{lem:KSescape} applies
to the restricted LP: at a feasible basis with BFS $x$ satisfying
$c\T x<c\T t$, the antistalling procedure with guide $t-x$ performs at
most $\min\{r(J)-1,m-1\}$ degenerate Simplex pivots followed by a
nondegenerate one.

\paragraph{The smallest face containing two points.}
The next lemma holds for every polytope in SEF.

\begin{lemma}\label{lem:commonzeros}
For every $x,y\in Q$,
$F(x,y)=\{z\in Q: z_i=0\ (i\in Z(x,y))\}$.
\end{lemma}
\begin{proof}
Let $Z=Z(x,y)$ and let $F_Z$ denote the right-hand side. The linear
function $\sum_{i\in Z}z_i$ is nonnegative on $Q$ and vanishes at $x$,
so $F_Z$, its set of minimizers over $Q$, is a face of $Q$ containing
$x$ and $y$; hence $F(x,y)\subseteq F_Z$. Conversely, let $G$ be a face
containing $x$ and $y$, let $p=\frac12(x+y)\in G$, and note that
$p_i=0$ for $i\in Z$ and $p_i>0$ for $i\notin Z$. For $z\in F_Z$ and
$\varepsilon>0$ small enough, the point $z'=p+\varepsilon(p-z)$ belongs to
$Q$: $Az'=b$, the coordinates of $z'$ indexed by $Z$ are zero, and the
other coordinates of $p$ are positive. Since
$p=\frac{\varepsilon}{1+\varepsilon}z+\frac{1}{1+\varepsilon}z'$ lies in the
relative interior of $[z,z']\subseteq Q$ and $G$ is a face containing
$p$, we get $z\in G$. Hence $F_Z\subseteq G$.
\end{proof}

The next lemma is the only place where the $0/1$ structure is used in
the proof of Theorem~\ref{thm:main}. 

\begin{lemma}\label{lem:binary}
Let $x,y,t$ be vertices of $Q$ with $F(x,t)=F(y,t)$. Then $x=y$.
\end{lemma}
\begin{proof}
Let $i\in V$ with $x_i=t_i$. If $t_i=0$, the inequality $z_i\geq0$ is
valid for $Q$ and tight at $x$ and $t$; if $t_i=1$, the inequality
$z_i\leq1$ is valid for $Q$, because $\mathcal P\subseteq[0,1]^d$, and tight at
$x$ and $t$. In both cases the face of $Q$ defined by the inequality
contains $x$ and $t$, hence $F(x,t)=F(y,t)$ and in particular $y$, so
$y_i=t_i$. By symmetry, $x_i=t_i$ if and only if $y_i=t_i$; and if
$x_i\neq t_i$, then $x_i=1-t_i=y_i$ since all values are binary. Hence
$x_V=y_V$, and $x=y$ because a vertex of $Q$ is determined by its
binary part.
\end{proof}

\paragraph{Lock or grow.}
Let $B$ be a feasible basis with BFS $x$, and consider a nondegenerate
Simplex pivot from $B$ with entering index $j$, direction $r=r^j$ and
step length $\theta>0$. Let $x'=x+\theta r$ be the vertex reached, and
let $L=\{i\in B: r_i<0,\ x_i/(-r_i)=\theta\}$ be the set of indices
attaining the minimum in the ratio test~\eqref{eq:ratio}, any of which
can be chosen as the leaving index.

\newpage
\begin{lemma}\label{lem:alternative}
With the above notation, exactly one of the following holds.
\begin{enumerate}[label=(\alph*)]
\item Some $g\in L$ satisfies $t_g=0$. Then $x_g>0$, and if $g$ is
chosen as the leaving index, $g$ is nonbasic after the pivot and
$x'_g=t_g=0$.
\item Every $g\in L$ satisfies $t_g>0$. Then, whichever leaving index
is chosen, $Z(x',t)\subseteq Z(x,t)$ and $F(x,t)\subsetneq F(x',t)$.
\end{enumerate}
\end{lemma}
\begin{proof}
(a) For $g\in L$ with $t_g=0$ we have $x_g=\theta(-r_g)>0$ and
$x'_g=x_g+\theta r_g=0$; if $g$ leaves, it is nonbasic after the pivot.

(b) Let $i\in Z(x',t)$, and suppose that $x_i>0$. Then $i\in B$, and
$x'_i=x_i+\theta r_i=0$ implies $r_i<0$ and $x_i/(-r_i)=\theta$, so
$i\in L$ with $t_i=0$, contradicting the assumption of (b). Hence
$x_i=0$ and $i\in Z(x,t)$, which proves $Z(x',t)\subseteq Z(x,t)$, and
Lemma~\ref{lem:commonzeros} yields $F(x,t)\subseteq F(x',t)$. Finally,
$x\neq x'$ because $x'_j=\theta>0=x_j$, and both are vertices of $Q$,
so $F(x,t)\neq F(x',t)$ by Lemma~\ref{lem:binary}.
\end{proof}

In case (a) the variable $g$ can be locked, which shrinks $Q_J$ while
keeping $t$ in it; in case (b) $Q_J$ is unchanged and $F(x,t)$ strictly
grows inside it, which by Lemma~\ref{lem:commonzeros} means that some
common zero of $x$ and $t$, possibly a basic variable of value zero,
becomes positive at $x'$. This is why only nonbasic common zeros are
locked.

\begin{algorithm}[t]
\caption{Target procedure}\label{alg:target}
\begin{algorithmic}[1]
\Require a feasible basis $B^0$ of~\eqref{eq:SEF} with BFS $x^0$, and
a vertex $t$ of $Q$ with $c\T t>c\T x^0$
\State $B\gets B^0$, $x\gets x^0$, $J\gets\varnothing$
\While{$c\T x<c\T t$}
\State run the antistalling procedure with guide $t-x$ on the LP
restricted to $[n]\setminus J$, stopping as soon as it selects an
entering index $j$ whose pivot is nondegenerate; let $L$ be the set of indices attaining the minimum ratio
in~\eqref{eq:ratio} for $j$
\If{some $g\in L$ satisfies $t_g=0$}
\State choose such a $g$ as the leaving index and set
$J\gets J\cup\{g\}$
\Else
\State choose any $g\in L$ as the leaving index
\EndIf
\State perform the pivot with entering index $j$ and leaving index
$g$, and update $B$ and $x$ accordingly
\EndWhile
\State \Return $B$
\end{algorithmic}
\end{algorithm}

We are now ready to prove Theorem~\ref{thm:target}. 

\begin{namedproof}{of Theorem~\ref{thm:target}}
The sequence of Simplex pivots is produced by
Algorithm~\ref{alg:target}, which repeatedly runs the antistalling
procedure on the restricted LP up to and including its first
nondegenerate pivot, choosing the leaving index of that pivot as
prescribed by Lemma~\ref{lem:alternative}. We call one execution of
the body of the while loop an \emph{iteration}.

Condition~\eqref{eq:J} holds initially and is maintained: $J$ only
grows by the leaving index $g$ of a nondegenerate pivot with $t_g=0$,
and by Lemma~\ref{lem:alternative}(a) $g$ is nonbasic after the pivot
with $x_g=t_g=0$. By Proposition~\ref{lem:restrict}, the current $B$ is
a feasible basis of the restricted LP with BFS $x$, and $t$ is a BFS
of it. As long as $c\T x<c\T t$, Lemma~\ref{lem:KSescape} applies to
the restricted LP, so an iteration performs at most
\begin{equation}\label{eq:itercost}
 \min\{r(J)-1,m-1\}+1=\min\{r(J),m\}
\end{equation} Simplex pivots, exactly one of which is nondegenerate. By
Proposition~\ref{lem:restrict} all these pivots are Simplex pivots
of~\eqref{eq:SEF}.

Call a \emph{phase} a maximal sequence of consecutive iterations during
which $J$ does not change, and its \emph{parameter} the value
$r=r(J)$ during the phase. Throughout a phase of parameter $r$, the
current vertex $x$ and $t$ belong to $Q_J$ and $x\neq t$, hence
$1\leq \dim F(x,t)\leq\dim Q_J\leq r$ (we abbreviate the dimension with $\dim$). An iteration whose nondegenerate
pivot falls into case (b) of Lemma~\ref{lem:alternative} strictly
increases the integer $\dim F(x,t)$, while degenerate pivots do not
change $x$; so at most $r-1$ iterations of the phase fall into case
(b), and an iteration falling into case (a) ends the phase. Hence a
phase of parameter $r$ consists of at most $r$ iterations and,
by~\eqref{eq:itercost}, performs at most $r\min\{r,m\}$ 
Simplex pivots, at most $r$ of which are nondegenerate.
Each phase, except possibly the last, ends by adding one index to $J$,
which decreases $r(J)$ by one. 
Summing over all values of $r$, the procedure
terminates after at most $\sum_{r=1}^{n-m}r\min\{r,m\}=\Cost{n}{m}$
 Simplex pivots, at most $\sum_{r=1}^{n-m}r=(n-m)(n-m+1)/2$ of
which are nondegenerate, at a basis whose BFS satisfies
$c\T x\geq c\T t$. Each pivot requires a basis update, a ratio test, a
guide update and the comparison of the values $t_g$, $g\in L$, with
zero, hence a polynomial number of arithmetic operations.
\end{namedproof}

Note that optimality of the target is not used: the proof only uses
that $t$ is a vertex with $c\T t>c\T x^0$, so that $t-x$ is an
improving feasible direction at every current vertex $x$ with
$c\T x<c\T t$.

We also remark that the sequence of Theorem~\ref{thm:target} is
produced by an explicit procedure, Algorithm~\ref{alg:target}, which,
given the target $t$, determines each pivot with a polynomial number
of arithmetic operations. For rational data, the numbers involved have
polynomial encoding length, so that each pivot is computed in
polynomial time: the BFS, the basic directions and the reduced costs
are obtained from nonsingular submatrices of $A$, and after every
degenerate pivot within one invocation of the antistalling procedure
each nonbasic coordinate of the guide is either zero or equal to
its value in $t-x$ (among the nonbasic coordinates before the pivot,
an update $y+\alpha r^f$ changes only $y_f$, which
becomes basic, and the leaving index has guide value zero, see the
proof of Theorem~1 in~\cite{KS}), which determines the basic
coordinates $y_B=-A_B^{-1}A_{\NN}y_{\NN}$ and $\alpha$.

\section{Optimal-basis cleanup}\label{sec:cleanup}
To prove our final bound, it only remains to discuss moving to an optimal basis once an optimal vertex has been reached.
In this context, Lemma~4 of~\cite{KS} states that, from a nonoptimal basis representing
an optimal vertex, a prescribed optimal basis $B^*$ of the same vertex can
be reached by at most $n-m$ degenerate Simplex pivots. The statement
however is slightly incorrect: reaching $B^*$ cannot in general be enforced (see Example 1 below).
Furthermore, after a degenerate pivot with leaving index $i$, the proof removes the variable
$x_i$ from the LP and continues with the smaller LP, but optimality of a basis of the smaller LP does
not imply optimality for the original LP.

However, it is still true that $n-m$ degenerate Simplex pivots suffice to reach \emph{some} optimal basis.
We here report a self-contained proof, showing that degenerate pivots can be obtained by a dual-feasibility
procedure, as in~\cite{FT}.

\begin{lemma}\label{lem:KScleanup}
Let $x^*$ be an optimal vertex of~\eqref{eq:LP} and let $B$ be a
feasible basis representing $x^*$. Then an optimal basis can be reached
from $B$ by at most $n-m$ degenerate Simplex pivots.
\end{lemma}
\begin{proof}
Let $I=\supp(x^*)$. A vector $p\in\R^m$ is an optimal solution of
the dual $\min\{b\T p: A\T p\geq c\}$ of~\eqref{eq:LP} if and only if
its slack vector $u(p)=A\T p-c$ satisfies $u(p)\geq\zero$ and
$u(p)_I=\zero$ (dual feasibility and complementary slackness with
$x^*$). For such $p$ let 
\[I(p)=\{i\in[n]: u(p)_i=0\}\supseteq I\] be the
set of tight dual constraints. For a basis $B$ let
$p^B=A_B^{-\mathsf T}c_B$. Then, $u(p^B)_B=\zero$ and
$u(p^B)_j=-\bar c_j$ for $j\in\NN$, so $B$ is optimal if and only if
$p^B$ is dual feasible.

We maintain a feasible basis $B$ representing $x^*$, so that
$I\subseteq B$, together with a dual optimal solution $p$, and we
measure progress by $\Phi=|B\setminus I(p)|$, the number of basic
indices whose dual constraint is not tight at $p$. Writing $u=u(p)$,
we have $c=A\T p-u$ and $Ar^j=\zero$, hence for $j\in\NN$
\begin{equation}\label{eq:rcdual}
 \bar c_j=c\T r^j=-u\T r^j=-u_j-\sum_{i\in B\setminus I(p)}u_i\,r^j_i ,
\end{equation}
because $r^j_j=1$, $r^j_i=0$ for $i\in\NN\setminus\{j\}$ and $u_i=0$
for $i\in I(p)$.

\smallskip
Consider the following:
\begin{enumerate}[label=(\roman*)]
\item If $B\subseteq I(p)$, then $A_B\T p=c_B$, so $p=p^B$ is dual
 feasible and $B$ is optimal.
\item If some $f\in I(p)\setminus B$ has $\bar c_f>0$, then $u_f=0$
 and~\eqref{eq:rcdual} give $r^f_i<0$ for some $i\in B\setminus I(p)$.
 As $i\notin I$, we have $x^*_i=0$, so $i$ attains the minimum ratio
 $0$ in~\eqref{eq:ratio}: entering $f$ and leaving $i$ is a degenerate
 Simplex pivot, the new basis represents $x^*$, and $\Phi$ decreases
 by one.
\item If $B$ is not optimal and no index as in (ii) exists, we move
 $p$ towards $p^B$. For $\tau\in[0,1]$ let $p(\tau)=p+\tau(p^B-p)$.
 Setting $\bar c_i=0$ for $i\in B$ gives slacks
 \[
 u(p(\tau))_i=(1-\tau)u_i-\tau\bar c_i.
 \] 
 For $i\in I(p)$, $u(p(\tau))_i$ equals
 $-\tau\bar c_i\geq0$, as $\bar c_i=0$ if $i\in B$ and $\bar c_i\leq0$
 otherwise by assumption. For $i\notin I(p)$ with $\bar c_i\leq0$, $u(p(\tau))_i$ is
 positive for $\tau<1$. The remaining indices have $\bar c_i>0$; they
 exist because $B$ is not optimal, they are nonbasic and not in
 $I(p)$, and their slack vanishes at $\tau_i=u_i/(u_i+\bar c_i)\in(0,1)$.
 Let $f$ attain $\tau'=\min\tau_i$ and let $p'=p(\tau')$. Then $p'$ is
 dual optimal, $I(p')\supseteq(B\cap I(p))\cup\{f\}$ and
 $I(p')\cap(B\setminus I(p))=\varnothing$. Hence replacing $p$ by $p'$
 does not change $\Phi$, and $f\in I(p')\setminus B$ has
 $\bar c_f>0$, so that (ii) applies.
\end{enumerate}
Start with $p=p^{B^*}$ for an optimal basis $B^*$ (one exists), which
is dual optimal with $I(p)\supseteq B^*$, so that
initially $\Phi\leq|B\setminus B^*|=|B^*\setminus B|\leq n-m$. Every
pivot decreases $\Phi$ by one and, by (i), $B$ is optimal once
$\Phi=0$; hence at most $n-m$ degenerate Simplex pivots are performed.
\end{proof}

\paragraph{Example 1.}
Consider $\max\{c\T x:Ax=b,\ x\geq\zero\}$ with
\[
 A=\begin{pmatrix}
 1&1&1&1\\
 0&1&-1&1
 \end{pmatrix},\qquad
 b=\zero,\qquad c=(0,0,-1,1)^{\mathsf T}.
\]
The matrix $A$ has rank $2$, and the first equality together with
nonnegativity implies $Q=\{\zero\}$. Hence $x^*=\zero$ is an optimal
vertex. Let $B=\{1,2\}$ and prescribe $B^*=\{1,3\}$; both are
feasible bases representing $x^*$.

At $B$, the nonbasic reduced costs are
$
 \bar c_3=-1,\mbox{ and } \bar c_4=1.
$
Thus $B$ is nonoptimal and only variable $4$ can enter in an improving
pivot. Since $A_4=A_2$, its basic direction is
$r^4=(0,-1,0,1)^{\mathsf T}$. The unique leaving variable is $2$, and
the step length is zero. The resulting basis is $B'=\{1,4\}$, with
nonbasic reduced costs
$
 \bar c_2=-1,\mbox{ and } \bar c_3=0.
$

Therefore $B'$ is optimal and no further improving pivot is possible.
The prescribed basis $B^*$ is also optimal, since its nonbasic reduced
costs are $\bar c_2=-1$ and $\bar c_4=0$. Nevertheless, $B^*$ is
unreachable from $B$ by Simplex pivots.

The example shows why the conclusion must allow the final optimal
basis to differ from the supplied guide.

\section{Final bound and simulation of monotone paths}\label{sec:algorithm}

Putting things together, we can now prove Theorem~\ref{thm:main}.

\begin{namedproof}{of Theorem~\ref{thm:main}}
Consider the standard equality form~\eqref{eq:SEF} of the $0/1$-LP,
with $n$ variables and $m$ constraints, so that $n-m=d-e\leq d$, and
let $t$ be an optimal vertex of $Q$. 

If the BFS $x^0$ of the given
basis $B^0$ satisfies $c\T x^0<c\T t$, Theorem~\ref{thm:target} with target
$t$ reaches a feasible basis $B$
whose BFS $x$ satisfies $c\T x\geq c\T t$, so that $x$ is an optimal
vertex, after at most $\Cost{n}{m}$ Simplex pivots. Among these, at most
$(n-m)(n-m+1)/2$ are nondegenerate.  
Otherwise, let $B=B^0$. By Lemma~\ref{lem:KScleanup}, an
optimal basis can be reached from $B$ by at most $n-m$ further
degenerate Simplex pivots. 

In total, at most
$\Cost{n}{m}+(n-m)\leq\sum_{r=1}^{d}r^2+d=O(d^3)$ Simplex pivots are
performed, at most $(n-m)(n-m+1)/2=O(d^2)$ of which are
nondegenerate.
\end{namedproof}

We now observe that Theorem~\ref{thm:target} can be applied repeatedly, with the vertices
of any short monotone path as successive targets. This simulates the
path in terms of objective values, with an overhead of
$\Cost{n}{m}=O(d^3)$ pivots per step, although the vertices actually
visited may be different.

\begin{corollary}\label{thm:simulate}
Suppose we are given a successor oracle which, queried at a vertex $v$
of $Q$, either returns a neighbor of $v$ with larger objective value
or reports that $v$ is optimal, and let $x^0,x^1,\ldots,x^L$ be the
path obtained by querying it repeatedly from $x^0$, so that $x^L$ is
optimal. 
From any feasible basis representing $x^0$, a basis whose BFS is
optimal can be reached by interleaving oracle calls and Simplex
pivots, so that for every $i$ the oracle is queried at $x^i$ only once
a BFS with objective value at least $c\T x^i$ has been reached, which
happens after at most $i\cdot\Cost{n}{m}$ Simplex pivots. An optimal
basis can then be reached with at most $n-m$ further degenerate
Simplex pivots.

\end{corollary}

\begin{proof}
Before querying $x^i$, maintain a BFS $w$ with $c\T w\geq c\T x^i$,
initially $w=x^0$. Upon receiving $x^{i+1}$, apply
Theorem~\ref{thm:target} with $x^{i+1}$ as target if $c\T w<c\T x^{i+1}$ (otherwise keep
the basis unchanged). This restores the invariant using at most
$\Cost{n}{m}$ additional pivots, so induction gives at most
$i\cdot\Cost{n}{m}$ pivots before querying $x^i$.
When the oracle reports that $x^L$ is optimal, $w$ is also optimal.
Lemma~\ref{lem:KScleanup} supplies the final at most $n-m$
degenerate Simplex pivots.
\end{proof}

Corollary~\ref{thm:simulate} reduces the construction of a short sequence
of Simplex pivots to that of a short monotone path. The path need not
be known in advance: when a local successor rule is available, the
target $x^{i+1}$ can be computed from $x^i$, for example using the constraints tight
there.

Using the above corollary, our result can be combined nicely with the $\ell_1$-steepest-edge rule and the shadow rules of~\cite{BDLKS,DLS,BKS}.
For several such rules (in particular, the classical $\ell_1$-steepest-edge rule) the next edge can be computed by optimizing a linear
function over an appropriate normalized section of the cone
of feasible directions at the current vertex~\cite{DLS,BDLKS}. For rational data, these local LPs can be solved in strongly polynomial time, exploiting the 0/1 structure through the method of~\cite[Section~6.6]{GLS}.
While one should note that all $0/1$-LPs with rational data can be solved in strongly polynomial time with the exact same technique,
 faster methods for solving the local LPs are often available in several combinatorial optimization
settings (e.g., matching and flows). Together with the remark
following Theorem~\ref{thm:target} and with Lemma~\ref{lem:KScleanup}
(whose proof is constructive once an optimal basis is known), 
these path rules thus yield pivot rules, with a polynomial number of
arithmetic operations per pivot.

\paragraph{AI Disclosure.} 
The proofs of all the results arose after multiple discussions with ChatGPT (GPT 6 Pro). Most notably, ChatGPT came up with constructions (not included in this paper) showing that, for prescribed bases, steepest edge is not approximable even with an irredundant description, and this was instrumental in moving to different ways of measuring progress. For editing, both ChatGPT and Claude have been used. 
The authors take full responsibility for the content and writing of the paper.


\newpage
\begin{thebibliography}{14}
\bibitem{BDLKS}
Black, A.E., De Loera, J.A., Kafer, S., Sanit\`a, L.: On the Simplex
method for $0/1$-polytopes. Mathematics of Operations Research
\textbf{50}(2), 1398--1420 (2025).
\url{https://doi.org/10.1287/moor.2021.0345}

\bibitem{BKS}
Black, A.E., Kafer, S., Sanit\`a, L.: Escaping degeneracy by following
shadow edges. arXiv:2610.12130 (2026).
\url{https://arxiv.org/abs/2610.12130}

\bibitem{Chubanov}
Chubanov, S.: A scaling algorithm for optimizing arbitrary functions
over vertices of polytopes. Mathematical Programming \textbf{190},
89--102 (2021). \url{https://doi.org/10.1007/s10107-020-01522-0}

\bibitem{Dantzig}
Dantzig, G.B.: Maximization of a linear function of variables subject
to linear inequalities. In: Koopmans, T.C. (ed.) Activity Analysis of
Production and Allocation, pp. 339--347. Wiley, New York (1951)

\bibitem{DLS}
De Loera, J.A., Kafer, S., Sanit\`a, L.: Pivot rules for
circuit-augmentation algorithms in linear optimization. SIAM Journal on
Optimization \textbf{32}(3), 2156--2179 (2022).
\url{https://doi.org/10.1137/21M1419994}

\bibitem{DS}
Dongarra, J., Sullivan, F.: Guest editors' introduction: The top 10
algorithms. Computing in Science \& Engineering \textbf{2}(1), 22--23
(2000). \url{https://doi.org/10.1109/MCISE.2000.814652}

\bibitem{FT}
Fukuda, K., Terlaky, T.: On the existence of a short admissible pivot
sequence for feasibility and linear optimization problems.
Pure Mathematics and Applications \textbf{10}(4), 431--447 (2000).
\url{https://www.cs.mcgill.ca/~fukuda/download/paper/short000214.pdf.gz}

\bibitem{GLS}
Gr\"otschel, M., Lov\'asz, L., Schrijver, A.: Geometric Algorithms and
Combinatorial Optimization, Algorithms and Combinatorics, vol.~2.
Springer, Berlin, 2nd edn. (1993).
\url{https://doi.org/10.1007/978-3-642-78240-4}

\bibitem{KM}
Kitahara, T., Mizuno, S.: A bound for the number of different basic
solutions generated by the Simplex method. Mathematical Programming
\textbf{137}(1--2), 579--586 (2013).
\url{https://doi.org/10.1007/s10107-011-0482-y}

\bibitem{KS}
Kukharenko, K., Sanit\`a, L.: On the number of degenerate Simplex
pivots. Mathematical Programming \textbf{216}, 757--791 (2026).
\url{https://doi.org/10.1007/s10107-026-02349-x}, see also
\url{https://arxiv.org/abs/2311.15799}

\bibitem{MatsuiTamura}
Matsui, T., Tamura, S.: Adjacency on combinatorial polyhedra. Discrete
Applied Mathematics \textbf{56}(2--3), 311--321 (1995)

\bibitem{Naddef}
Naddef, D.: The Hirsch conjecture is true for $(0,1)$-polytopes.
Mathematical Programming \textbf{45}, 109--110 (1989).
\url{https://doi.org/10.1007/BF01589099}

\bibitem{Santos}
Santos, F.: A counterexample to the Hirsch conjecture. Annals of
Mathematics \textbf{176}(1), 383--412 (2012).
\url{https://doi.org/10.4007/annals.2012.176.1.7}

\bibitem{Smale}
Smale, S.: Mathematical problems for the next century. In: Mathematics:
Frontiers and Perspectives. American Mathematical Society (2000)
\end{thebibliography}
\end{document}